Two Sumerian Words of Fractions in Babylonian Mathematics: igi-n-gál and igi-te-en

Kazuo MUROI

§1. Introduction

In Babylonian mathematics two Sumerian words of fractions occur, which were originally used in non-mathematical texts. They are igi-n-gál "the reciprocal of (the number) n", which is often abbreviated to igi-n, and igi-te-en whose meaning is somewhat abstract, that is, "a proportion (of something)" or "the ratio (of something to another)"[1]. Thus the mathematical meanings of the two terms are quite clear, but we have not been able to clarify their literal meanings or their origins so far. A major cause of the difficulty in studying the etymologies of the terms would lie in the following facts:

In spite of numerous administrative and economic documents from the 25th century B.C. on, where igi-n-gál occurs, we cannot find any instance among them that contributes to the clarification of the original meaning of the term. As for igi-te-en, we have only a small number of instances of the term in legal texts and in literary texts of the Old Babylonian period (ca. 2000-1600 B.C.) which also provide no clue to our problem.[2]

Therefore, the only thing left for us to do is search mathematical texts of the Old Babylonian period for the keys to the enigmatic origins of igi-n-gál and igi-te-en.

In the present paper I shall offer a most probable interpretation of the early term, igi-n-gál, and a definite solution to the etymology of the later term, igi-te-en, both of which are based on my analysis of several mathematical terms that concern multiplication or division.

§2. The origin of igi-n-gál

In many documents of around the 25th century B.C. the unit fractions igi-3-gál (= 1/3), igi-4-gál (= 1/4), and igi-6-gál (= 1/6) as well as the so-called natural fractions bar (= 1/2), šušana (= 1/3), and šanabi (= 2/3) are frequently used concerning the weight of silver. I cite a legal document of the 24th century B.C. as an example:

3 1/3$^{ša}$ (= šušana) ma-na 1 gín igi-3-gál kù luh-ha / šeš-tur-ré / u$_4$-lú / e-da-tuku / u$_4$-ba / ur-é-mùš-ke$_4$ / di-bi ì-ku$_5$ / lugal-an-da / ensí-kam

"Šeštur lent 3 1/3 mana 1 1/3 shekel of refined silver to Ulu. At that time Uremuš passed judgment on it. Lugalanda was the city-ruler (of Lagash)."[3]

Since the meaning of the noun igi is "eye, face, front" and that of the verb gál is "to be, to place, to have", the most probable literal meaning of igi-n-gál would be "that which exists in front of n", which will be also supported by mathematical evidence.

In Babylonian mathematics too the multiplication of integers is defined by repeated addition as the following typical expression shows:

m *a-na* n *našûm* or m a-rá n túm "to carry m n times", that is, "to multiply m by n".[4]

This was derived from an everyday expression of those days such as "to carry (something) n times"

or "to carry (something) repeatedly", and it should be noted that the multiplicand m usually precedes the multiplier n in its word order. We also find the same word order of multiplicand and multiplier in some mathematical texts where the division of a number by an irregular number is performed, for example:

igi-40,51 *ú-la ip-pa-ṭa-ar mi-nam a-na* 40,51 *lu-uš-ku-un ša* 10,25,45 *i-na-di-nam* 15 ba-an-da-*šu*

"The reciprocal of 40,51 (= 3·19·43) is not obtained. What should I put down for 40,51 so that it gives me 10,12;45? 0;15, (by) its discernment."[5]

If we multiply 0;15 by 40,51 after the Babylonians:

0;15 *a-na* 40,51 *našûm* "to multiply 0;15 by 40,51",

we will obtain the product 10,12;45 and eventually solve the linear equation 40,51x = 10,12;45. Note that the interrogative *mīnam* (the accusative of *mīnum* "what") or the number 0;15 precedes the multiplier 40,51.

Therefore it is reasonable to assume that the Sumerians of third millennium B.C. defined one third, for example, as "that which makes one if it is carried three times". In other words, they considered a relation, 1/3·3=1, where 1/3 is the multiplicand and 3 the multiplier as in the case of the Japanese language in contrast to English. They must have named one third igi-3-gál "that which exists in front of 3" besides šušana. This interpretation of igi-n-gál is also supported by its Akkadian counterpart, *pāni*-n "the reciprocal of n, (literally) the front of n" which must be a literal translation of the abbreviated Sumerian term igi-n.

§3. The origin of igi-te-en

Although the noun igi of igi-n-gál is "front", as has been shown in §2, the igi of igi-te-en is, in my judgment, not "front" but "bubble". Through some lexical lists we know the fact that the cuneiform sign IGI occasionally means "foam, bubble" which is derived from its usual meaning "eye":

im-hu-ur IGI-A = *hu-ur-hu-ma-at me-e* "foam on water" (IGI-A = imhur$_4$)
im-hu-ur IGI-KAŠ = *hu-ur-hu-ma-at ši-ka-ri-im* "foam on beer" (IGI-KAŠ = imhúr)
im-hu-ur IGI-GA = *hu-ur-hu-mat* ga-meš "foam on milk" (IGI-GA = imhùr).[6]

As to the verb te-en or te(n) "to extinguish, to cool down" we have no problem, and so we may translate igi-te-en as "to deflate a foam" which will be also mathematically supported in the next section. Moreover, the Akkadian verb *belûm* "to become extinguished, to come to an end, to burst (said of bubbles)"[7], which is one of the Akkadian counterparts of the Sumerian verb te-en, occurs concerning a burst of bubble in Old Babylonian oil omen texts, for example:

*šumma* (DIŠ) *šulman* (SILIM) *id-di-a-am-ma ù be-li mar-ṣum i-ma-at ṣāb* (ERÍN) *harrānim* (KASKAL) *re-eš eqli* (A-ŠÀ)-*šu ú-ul i-ka-aš-ša-ad*

"If (the oil) produces a bubble and it bursts, the sick person will die, (and) the army on an expedition will not reach its destination."[8]

In order to show that the Babylonians metaphorically described a broken piece or a fraction that resulted from some division as igi-te-en "deflated foam", I cite two passages from a famous literary work, *Lamentation over Sumer and Ur*[9]:

$^d$Nanna ùg u$_8$-gin$_7$ lu-a-na igi-te-en-bi si-il-le-dè
"to break into pieces Nanna's people numerous as sheep"
kur-kur-re du$_{10}$-ús dili dab$_5$-ba-bi igi-te-en-bi [ba-si-il]
"When all the countries held a single mode of life, they were broken into pieces."

If we carefully observe that a bubble bursts and often breaks into smaller bubbles when we are in cooking, washing, or drinking beer, we could more clearly understand the original meaning of the igi-te-en probably coined in the early years of the Old Babylonian period.

§4. The antonym of a-rá-kár

In Babylonian mathematics the successive doubling of a number, which was especially called a-rá-kár "blowing up multiplication" in Sumerian, was used together with the successive halving of its reciprocal in order to make large tables of reciprocals. If we look at one of the tables, for example, a table of the reciprocal pairs of and (n=1 2 3 … 30), we will immediately see that how the numbers are blowing up in succession and how the reciprocals are diminishing in succession. See Table 1. There is a possibility that the latter successive halving of a reciprocal was called igi-te-en "deflated foam" in connection with the a-rá-kár. We may have a piece of evidence in a Sumerian proverb that the verbs bar$_7$ (≈ kár) and te-en occur side by side. The verb bar$_7$ as well as the verb kár corresponds to the Akkadian verb *napāhum* "to blow (something), to light a fire, (in the stative) to be bloated, swollen, blown up"[10]:

hé-lu-a nam-ba-lá / hé-diri na-an-su-su / hé-ib-bar$_7$ na-an-te-en-te-en

"Let it be abundant, do not diminish it.  Let it be surplus, do not restore it.  Let it be swollen, do not make it go down."[11]

In conclusion, the antonym of the mathematical term a-rá-kár is probably our igi-te-en discussed in §3.

Notes

Table

| n | $2^n \cdot 5$ | |
|---|---|---|
| 1 | 10 | 0;6 |
| 2 | 20 | 0;3 |
| 3 | 40 | 0;1,30 |
| 4 | 1,20 | 0;0,45 |
| 5 | 2,40 | 0;0,22,30 |
| 6 | 5,20 | 0;0,11,15 |
| 7 | 10,40 | 0;0,5,37,30 |
| 8 | 21,20 | 0;0,2,48,45 |
| 9 | 42,40 | 0;0,1,24,22,30 |
| 10 | 1,25,20 | 0;0,0,42,11,15 |
| 11 | 2,50,40 | 0;0,0,21,5,37,30 |
| 12 | 5,41,20 | 0;0,0,10,32,48,45 |
| 13 | 11,22,40 | 0;0,0,5,16,24,22,30 |
| 14 | 22,45,20 | 0;0,0,2,38,12,11,15 |
| 15 | 45,30,40 | 0;0,0,1,19,6,5,37,30 |
| 16 | 1,31,1,20 | 0;0,0,0,39,33,2,48,45 |
| 17 | 3,2,2,40 | 0;0,0,0,19,46,31,24,22,30 |
| 18 | 6,4,5,20 | 0;0,0,0,9,53,15,42,11,15 |
| 19 | 12,8,10,40 | 0;0,0,0,4,56,37,51,5,37,30 |
| 20 | 24,16,21,20 | 0;0,0,0,2,28,18,55,32,48,45 |
| 21 | 48,32,42,40 | 0;0,0,0,1,14,9,27,46,24,22,30 |
| 22 | 1,37,5,25,20 | 0;0,0,0,0,37,4,43,53,12,11,15 |
| 23 | 3,14,10,50,40 | 0;0,0,0,0,18,32,21,56,36,5,37,30 |
| 24 | 6,28,21,41,20 | 0;0,0,0,0,9,16,10,58,18,2,48,45 |
| 25 | 12,56,43,22,40 | 0;0,0,0,0,4,38,5,29,9,1,24,22,30 |
| 26 | 25,53,26,45,20 | 0;0,0,0,0,2,19,2,44,34,30,42,11,15 |
| 27 | 51,46,53,30,40 | 0;0,0,0,0,1,9,31,22,17,15,21,5,37,30 |
| 28 | 1,43,33,47,1,20 | 0;0,0,0,0,0,34,45,41,8,37,40,32,48,45 |
| 29 | 3,27,7,34,2,40 | 0;0,0,0,0,0,17,22,50,34,18,50,16,24,22,30 |
| 30 | 6,54,15,8,5,20 | 0;0,0,0,0,0,8,41,25,17,9,25,8,12,11,15 |